\documentclass{article}
\usepackage{amssymb}
\input diagrams
\title{Combinatorial Structure of Finite \\ 
Dimensional Representations of Yangians: \\
the Simply-Laced Case}  % yuk
\author{Michael Kleber\thanks{Supported by NSF grant DMS 94-01163.}}
\date{}
\newtheorem{thm}{Theorem}
\newtheorem{cor}[thm]{Corollary}
\newtheorem{lemma}[thm]{Lemma}

\renewcommand{\d}{{\mathbf d}}
\newcommand{\dsum}{\oplus}
\newcommand{\Dsum}{\bigoplus}
\newcommand{\Dt}{{\Delta_1\ldots\Delta_t}}
\newcommand{\Dtmt}{{\Delta_1^{m_1}\!\!\ldots\Delta_t^{m_t}}}
\newcommand{\g}{{\mathfrak g}}
\newcommand{\ghat}{{\hat{\mathfrak g}}}
\newcommand{\hyt}{{\mathop{\mathrm{ht}}\nolimits}}
\newcommand{\iso}{\simeq}
\renewcommand{\l}{\ell}
\renewcommand{\O}{{\mathcal O}}
\newcommand{\oplustag}[1]{\mathbin{\mathop{\oplus}\limits^{#1}}}
\newcommand{\R}{{\mathbb R}}

\newcommand{\tl}{{T(\l)}}

\newcommand{\wl}{{\omega_\l}}
\newcommand{\wml}{{W_m(\l)}}

\begin{document}
\maketitle

\begin{abstract}
We compute the decomposition of representations of Yangians into
$\g$-modules for simply-laced $\g$.  The decomposition has an interesting
combinatorial tree structure.  Results depend on a conjecture of Kirillov
and Reshetikhin.
\end{abstract}
\addtocounter{section}{-1} % so introduction is section 0
\section{Introduction}
\label{sec_intro}

Let $\g$ be a complex semisimple Lie algebra of rank $r$, and $Y(\g)$ its
Yangian (\cite{Dr}), a Hopf algebra which contains the universal
enveloping algebra $U(\g)$ of $\g$ as a Hopf subalgebra.  Write
$\alpha_1,\ldots,\alpha_r$ for the fundamental roots and
$\omega_1,\ldots,\omega_r$ for the fundamental weights of $\g$.  As
defined in~\cite{KR}, denote by $\wml$ a particular irreducible $Y(\g)$
module all of whose $\g$-weights $\lambda$ satisfy $\lambda\preceq m\wl$,
where $\alpha\preceq\beta$ means $\beta-\alpha$ is a nonnegative integer
linear combination of the roots $\{\alpha_i\}$.  Specifically, $\wml$
decomposes into $\g$-modules as
\begin{equation}
\label{def_decomp}
\wml |_{\g} \iso \Dsum_{\lambda\preceq m\wl} V_\lambda^{\dsum n_\lambda}
\end{equation}
where $V_\lambda$ is the irreducible $\g$-module with highest weight
$\lambda$ and it occurs $n_\lambda$ times in $\wml$.  In particular,
$n_{m\wl}=1$.

There is a formula for the multiplicities $n_\lambda$ in~\cite{KR} based
on the conjecture that every finite-dimensional representation of $Y(\g)$
can be obtained from one specific representation by means of the
``reproduction scheme,'' defined in~\cite{KRS}.  If this conjecture holds,
then write $\lambda = m\wl - \sum n_i \alpha_i$, and it is proved
in~\cite{KR} that
\begin{equation}
\label{defZ}
n_\lambda = Z(\l,m|n_1,\ldots,n_r) =
\sum_{\mbox{partitions}} \;\; \prod_{n\geq1} \;\; \prod_{k=1}^r
{{P^{(k)}_n(\nu) + \nu^{(k)}_n} \choose {\nu^{(k)}_n}}
\end{equation}
The sum is taken over all ways of choosing partitions
$\nu^{(1)},\ldots,\nu^{(r)}$ such that $\nu^{(i)}$ is a partition of
$n_i$ which has $\nu^{(i)}_n$ parts of size $n$ (so $n_i =
\sum_{n\geq1} n \nu^{(i)}_n$).  The function $P$ is defined by
\begin{eqnarray}
\label{defPgen}
P^{(k)}_n(\nu) &=& \min(n,m)\delta_{k,\l}
  - 2 \sum_{h\geq 1} \min(n,h)\nu^{(k)}_{h} + \\
&&\hspace{1cm} +
  \sum_{j\neq k}^r \sum_{h\geq 1} \min(-c_{k,j}n,-c_{j,k}h)\nu^{(j)}_{h}
\nonumber
\end{eqnarray}
where $C=(c_{i,j})$ is the Cartan matrix of $\g$.  We define
${a\choose b}$ to be 0 whenever $a<b$; since the values of $P$ can be
negative, many of the binomial coefficients in~(\ref{defZ}) can be zero.

Yangians are closely related to quantum affine universal enveloping
algebras $U_q(\ghat)$ when $q$ is not a root of unity, and $U_q(\g)$ is a
Hopf subalgebra of $U_q(\ghat)$ in much the same way that $\g$ is a Hopf
subalgebra of $Y(\g)$.  View $V_\lambda$ as a highest weight module over
$U_q(\g)$; then an affinization of $V_\lambda$ is defined in~\cite{ChP}
to be a $U_q(\ghat)$-module all of whose weights as a $U_q(\g)$-module
satisfy $\mu\preceq\lambda$ and $\lambda$ appears with multiplicity 1.
In the case that $\lambda=m\wl$, $V_{m\wl}$ has a unique minimal
affinization (\cite{ChP}) with respect to a partial ordering defined in
\cite{Ch}, and it is believed~(\cite{Rtalk}) that the decomposition of
this minimal affinization into $U_q(\g)$-modules is the same as the
decomposition for the Yangian module $\wml$.

In Section~\ref{sec_algorithm}, we view the values of $P^{(k)}_n$ as the
coordinates of certain strings of weights of $\g$ which lie inside the
Weyl chamber.  This interpretation allows us to compute the values of
$n_\lambda$ much more efficiently.  Furthermore, the ``initial
substring'' relation on the labelling by strings of weights imposes the
structure of a rooted tree on the set of $\g$-modules which make up
$\wml$, rooted at $V_{m\wl}$ and with the children of any $V_\lambda$
having highest weights $\mu=\lambda-\delta$ with $\delta$ in the positive
root lattice.

In Section~\ref{sec_growth}, we use this added structure to study the
asymptotics of the dimension of $\wml$ as $m$ gets large, based on the
fact that the tree structure of $\wml$ lifts to $W_{m+1}(\l)$.  We show
that the conjecture implies that the dimension grows asymptotically to a
polynomial in $m$, and compute the degree of this polynomial for every
simply-laced $\g$ and choice of $\wl$.

In Section~\ref{sec_table} we give a list of the decompositions of $\wml$
for all simply-laced $\g$ and small values of $m$ as derived 
numerically from the conjecture, using the results of
Section~\ref{sec_algorithm}.  For any choice of $\g$, representations
$W_1(\l)$ are called fundamental representations, since every
finite-dimensional representation of $Y(\g)$ appears as a quotient of a
submodule of a tensor product of fundamental representations.  In the
context of $U_q(\ghat)$-modules, the decompositions of most of the
fundamental representations were calculated in~\cite{ChP} using
completely different techniques, and those calculations agree with ours.

A similar idea can be used to give a combinatorial interpretation to the
values in equations~(\ref{defZ}) and~(\ref{defPgen}) when $\g$ is not
simply-laced.  The resulting structure is not as regular as in the
simply-laced case, but should yeild similar results.

The author is grateful to N. Yu. Reshetikhin for suggestion of the
problem, discussions, support and encouragement during the development
and preparation of this paper.

\section{Structure in the simply-laced case}
\label{sec_algorithm}

Assume that our Lie algebra $\g$ of rank $r$ is simply-laced.  Then
equation~(\ref{defPgen}) becomes
\begin{equation}
\label{defP}
P^{(k)}_n(\nu) = \min(n,m)\delta_{k,\l} - \sum_{j=1}^r c_{j,k}
     \left( \sum_{h\geq 1} \min(n,h)\nu^{(j)}_{h} \right)
\end{equation}
Fix a highest weight $m\wl$, and pick an arbitrary $\nu = 
(\nu^{(1)},\ldots,\nu^{(r)})$, where each $\nu^{(i)}$ is a partition 
of some integer $n_i$.  Then for any nonnegative integer $n$, the 
values $(P^{(1)}_n,\ldots,P^{(r)}_n)$ can be thought of as the 
$\omega$-coordinates of some weight; define
$$
\mu_n = \sum_{k=1}^r P^{(k)}_n \omega_k
$$
A given $\nu$ contributes a nonzero term to the sum in~(\ref{defZ})
if and only if the corresponding weights $\mu_0=0, \mu_1, \mu_2,\ldots$ 
all lie in the dominant Weyl chamber.
The motivation for seeing these as weights is that the sum
in~(\ref{defP}) can be naturally realized as subtracting some linear
combination of roots; if we let
\begin{equation}
\label{defd}
d_n = \sum_{k=1}^r
  \left( \sum_{h\geq 1} \min(n,h)\nu^{(k)}_{h} \right) \alpha_k
\end{equation}
then $\mu_n = \min(n,m)\wl - d_n$.

Think of $\nu^{(1)},\ldots,\nu^{(r)}$ as Young diagrams with $\nu^{(k)}$
having $\nu^{(k)}_{h}$ rows of length $h$.  Then we can tell whether a
sequence of vectors $d_0=0, d_1, d_2,\ldots$ can arise from
$\nu^{(1)},\ldots,\nu^{(r)}$ by looking at their successive differences
$\delta_i = d_i - d_{i-1}$.  If we write $\delta_n$ out as a linear
combination of the roots $\{\alpha_i\}$, then the $\alpha_k$-coordinate
is the number of boxes in the $n$th column of the Young diagram of
$\nu^{(k)}$, since the sum $\sum_h \min(n,h)\nu^{(k)}_{h}$
in~(\ref{defd}) is the number of boxes in the first $n$ columns.  Thus a
sequence arises from partitions if and only if the $\delta_i$ are
nonincreasing; that is, $\forall i\geq 1 : \delta_i \succeq \delta_{i+1}$.

If we let $s$ be the size of the largest part in any of the partitions 
in $\nu$, then $d_s = d_t$ for all $t>s$ (and $s$ is the smallest 
index for which this is true), and all the information we need to 
identify a particular summand of $\wml$ is the (strictly increasing) 
chain of weights $d_0=0 \prec d_1\prec\cdots\prec d_s$, which we 
define to have length $s$.  Note that the chain of length 0 consisting 
of only $d_0=0$ is permissible, arises from empty partitions, and 
corresponds to the $V_{m\wl}$ component of $\wml$.

In summary, we have proven the following:

\begin{thm}
\label{thm_decomp}
Let $\g$ be a simply-laced complex semisimple Lie algebra of rank $r$ 
with fundamental roots $\alpha_1,\ldots,\alpha_r$ and fundamental 
weights $\omega_1,\ldots,\omega_r$, and assume the decomposition of 
$\wml$ into $\g$-modules in (\ref{def_decomp}) is given by the 
conjecture in (\ref{defZ}) and (\ref{defPgen}).  Then that 
decomposition can be refined into a direct sum of parts indexed by 
chains of weights $\d=d_0,\ldots,d_s$ with successive differences 
$\delta_i = d_i - d_{i-1}$ (and $\delta_{s+1}=0$) such that
  \begin{enumerate}
  \item[(i)] $d_0=0$ and $d_0\prec d_1\prec\cdots\prec d_s$,
  \item[(ii)] $\min(n,m)\wl - d_n$ lies in the positive Weyl chamber 
  for $0\leq n\leq s$, and
  \item[(iii)] $\delta_i \succeq \delta_{i+1}$ for all $1\leq i\leq s$.
  \end{enumerate}
The summand with label $\d=d_0,\ldots,d_s$ consists of the $\g$-module 
of highest weight $m\wl - d_s$ with multiplicity
$$
\prod_{n\geq1} \;\; \prod_{k=1}^r \;
{{P^{(k)}_n(\d) + \d^{(k)}_n} \choose {\d^{(k)}_n}}
$$
where the values of $P^{(k)}_n(\d)$ and $\d^{(k)}_n$ are defined by the 
relations
\begin{eqnarray*}
\min(m,n)\wl - d_n &=& \sum_{k=1}^r P^{(k)}_n(\d) \omega_k \\
\delta_n - \delta_{n+1} &=& \sum_{k=1}^r \d^{(k)}_n \alpha_k
\end{eqnarray*}
and all of the multiplicities are nonzero.
\end{thm}

\begin{figure}
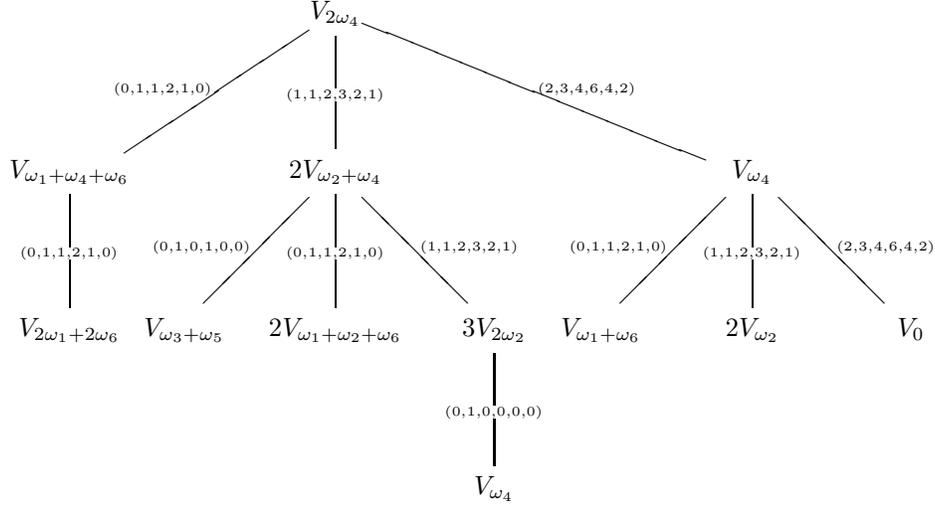

\begin{diagram}[labelstyle=\scriptscriptstyle,noPS]
&&&V_{2\omega_4}&&&&&&& \\
&& \ldLine(3,2) ^{(0,1,1,2,1,0)} & \dLine ~{(1,1,2,3,2,1)} & 
    \rdLine(5,2) ^{(2,3,4,6,4,2)} &&&&&& \\
V_{\omega_1 + \omega_4 + \omega_6}  & & &
    2 V_{\omega_2 + \omega_4} & & & & & V_{\omega_4} & & \\
\dLine ~{(0,1,1,2,1,0)} & & \ldLine ^{(0,1,0,1,0,0)} &
    \dLine ~{(0,1,1,2,1,0)} & \rdLine ^{(1,1,2,3,2,1)} & & &
    \ldLine ^{(0,1,1,2,1,0)} & \dLine ~{(1,1,2,3,2,1)} &
    \rdLine ^{(2,3,4,6,4,2)} & \\
V_{2\omega_1 + 2\omega_6} & \;\; V_{\omega_3 + \omega_5} & &
    2 V_{\omega_1 + \omega_2 + \omega_6} & & 3 V_{2\omega_2} &
    \;\; V_{\omega_1 + \omega_6} & & 2 V_{\omega_2} & & V_{0} \\
&&&&& \dLine ~{(0,1,0,0,0,0)} &&&&&\\
&&&&& V_{\omega_4} &&&&&\\
\end{diagram}
\caption{Tree-structure\protect\footnotemark\  of the 
decomposition of $W_2(4)$ for $E_6$.}
\label{fig_tree}
\end{figure}

This decomposition is a refinement of the one in~(\ref{def_decomp}) 
since it is possible to find two different chains $d_0,\ldots,d_s$ and 
$d'_0,\ldots,d'_t$ with $d_s = d'_t$.  This happens any time the sum 
in~(\ref{defP}) has more than one nonzero term.  One example of this 
occurs in $W_2(4)$ for $E_6$; see Figure~\ref{fig_tree}.

\begin{cor}
\label{cor_algor}
If $d_0,\ldots,d_s$ is a valid label then any initial segment 
$d_0,\ldots,d_{s'}$ (for $0 \leq s' < s$) is a valid label also.  
Conversely, given any label $d_0,\ldots,d_s$, we can extend it to 
another valid label by appending any weight $d_{s+1}$ which satisfies 
the conditions that $\min(s+1,m)\wl - d_{s+1}$ is in the 
positive Weyl chamber and, if $s>0$, that $d_s \prec d_{s+1} \preceq 
d_s+\delta_s$.
\end{cor}

This follows immediately from conditions {\em (i)--(iii)}.  Since $d_0$
must be 0, this completely describes an effective algorithm for computing
the conjectured decomposition of a given $\wml$.  The computations in
Section~\ref{sec_table} were computed using this algorithm.  The fact
that an initial segment of a valid label is still a valid label is the
key result which fails to hold true when $\g$ is not simply-laced.

Since truncating any label gives you another label, we can impose a 
tree structure on the parts of this decomposition, with a node of the 
tree corresponding to a summand in the decomposition from 
Theorem~\ref{thm_decomp}.  The ``children'' of the node with label 
$d_0,\ldots,d_s$ are all the nodes indicated by 
Corollary~\ref{cor_algor}; we can label the edges joining them to 
their parent with the various choices for the increment 
$\delta_{s+1}$.  For each $n\geq 0$, the $n$th row of the tree 
consists of all the nodes with labels of length $n$.

\footnotetext{Thanks to Paul Taylor's {\tt diagrams} package.}

%
% Warning: there's no guarantee that the footnote will end up on the same
% page as the figure!  I can't think of a way to avoid this problem.
%

As an example of this structure, the tree for the decomposition of 
$W_2(4)$ for $\g=E_6$ is given in Figure~\ref{fig_tree}.  Scalars in 
front of modules, as in $2 V_{\omega_2 + \omega_4}$, indicate 
multiplicity.  The label $(a_1,\ldots,a_6)$ corresponds to an 
increment $\delta=\sum a_i\alpha_i$, so condition {\em (iii)} says 
that the labels along any path down from $V_{2\omega_4}$ will be 
nonincreasing in each coordinate.  The labels on the edges are 
technically unnecessary, since they can be obtained by subtracting the 
highest weight of the child from the highest weight of the parent.  
However, as the next corollary shows, they do record useful 
information that is not apparent by looking directly at the highest 
weights.

\begin{cor}
\label{cor_lift}
If $d_0,\ldots,d_s$ is a valid label for $W_m(\l)$, then it is also a
valid label for $W_{m'}(\l)$ for any $m'>m$, and for any $m'\geq s$.
\end{cor}

Both parts are based on the fact that condition {\em (ii)} is the only 
one that depends on $m$.  For $m'>m$, if $\min(n,m)\wl - d_n$ is 
a nonnegative linear combination of the $\{\wl\}$ then adding 
some nonnegative multiple of $\wl$ will not change that fact.  
And if $m'\geq s$ then the value of $m'$ is irrelevant; the weights 
we look at are just $n\wl - d_n$ for $0\leq n\leq s$.

If we can lift labels from $W_m(\l)$ to $W_{m+1}(\l)$, we can also lift
the entire tree structure.  Specifically, the lifting of labels extends
to a map from the tree of $W_m(\l)$ to the tree of $W_{m+1}(\l)$ which
preserves the increment $\delta$ of each edge and lifts each $V_\lambda$
to $V_{\lambda+\wl}$.  The $m'\geq s$ part of 
Corollary~\ref{cor_lift} tells us that this map is a bijection on rows
$0,1,\ldots,m$ of the trees, where the labels have length $s\leq m$.  On
this part of the tree, multiplicities are also preserved.  This 
follows from the formula for multiplicities in Theorem~\ref{thm_decomp}:
the only values of $P^{(k)}_n(\d)$ that change are for $n=m+1$, but
$\d^{(k)}_n=0$ when $n$ is greater than the length of the label, so the
product of binomial coefficients is unchanged.

\section{The Growth of Trees}
\label{sec_tree}\label{sec_growth}

In this section, we will prove that as $m$ gets large, the dimension of
the representation $\wml$ grows like a polynomial in $m$, and will give a
method to compute the degree of the polynomial growth.  All statements 
assume the conjectural formulas for multiplicities of $\g$-modules.
Roots and weights are numbered as in~\cite{Bour}.

Since the tree decompositions for $\wml$ for $m=1,2,3,\ldots$ stabilize,
we can define $\tl$ to be the tree whose top $n$ rows coincide with those
of $\wml$ for all $m\geq n$.  The highest weight associated with an
individual node appearing in $\tl$ is only well-defined up to addition
of any multiple of $\wl$, but the difference $\delta$ between any
node and its parent is well-defined.  (These differences are the labels
on the edges of the tree in Figure~\ref{fig_tree}.)  We can characterize
each node by the string of successive differences $\delta_1
\succeq \delta_2 \succeq\cdots\succeq \delta_s$ which label the $s$ edges
in the path from the root of the tree to that node.  The multiplicity of
a node of $\tl$ is well-defined, as already noted.

The tree of $\wml$ matches $\tl$ exactly in the top $m$ rows.  The number
of rows in the tree of $\wml$ is bounded by the largest
$\alpha$-coordinate of $m\wl$, since if $\delta_1,\ldots,\delta_s$
is a label of $\wml$ then $m\wl - \sum_{i=1}^s \delta_i$ must be
in the positive Weyl chamber, and whatever $\alpha$-coordinate is nonzero
in $\delta_s$ must be nonzero in all of the $\delta_i$.  Therefore to
prove that the dimension of $\wml$ grows as a polynomial in $m$, it
suffices to prove that the dimension of the part of $\wml$ which
corresponds to the top $m$ rows of $\tl$ does so.

% Deleted paragraph used to be here.

Now we need to examine the structure of the tree $\tl$.  The path 
$\delta_1,\ldots,\delta_s$ to reach a vertex is a sequence of weights 
whose $\alpha$-coordinates are nonincreasing.  Write this instead as 
$\Dtmt$ where the $\Delta_i$ are strictly decreasing and $m_i$ is the 
number of times $\Delta_i$ occurs among $\delta_1,\ldots,\delta_s$; we 
will say this path has {\em path-type} $\Dt$.  The number of 
path-types that can possibly appear in the tree $\tl$ is finite, since 
each $\Delta_i$ is between $\wl$ and 0 and has integer 
$\alpha$-coordinates.

We need to understand which path-types $\Dt$ and which choices of
exponents $m_i$ correspond to paths which actually appear in $\tl$.
Given a path $\delta_1,\ldots,\delta_s$, assume that $m>s$ and recall
$\mu_n = n\wl - d_n = n\wl - \sum_{i=1}^n \delta_i$.
Condition~{\em (ii)} from Theorem~\ref{thm_decomp} requires that $\mu_n$
is in the positive Weyl chamber for $1\leq n\leq s$; that is, the
$\omega$-coordinates of $\mu_n$ must always be nonnegative.  (These
coordinates are just the values of $P^{(k)}_n$ from
Theorem~\ref{thm_decomp}.)  Since $\mu_n = \mu_{n-1} + \wl -
\delta_n$, we need to keep track of which $\omega$-coordinates of
$\wl-\delta_n$ are positive and which are negative.

For a path-type $\Dt$, we say that $\Delta_i$ {\em provides} 
$\omega_k$ if the $\omega_k$-coordinate of $\wl - \Delta_i$ is 
positive, and that it {\em requires} $\omega_k$ if the coordinate is 
negative.  Geometrically, $\Delta_i$ providing $\omega_k$ means that 
each $\Delta_i$ in the path moves the sequence of $\mu$s away from the 
$\omega_k$-wall of the Weyl chamber, while requiring $\omega_k$ moves 
towards that wall.  The terminology is justified by restating what 
condition {\em (ii)} implies about path-types in these terms:

\begin{lemma}
\label{lem_pathtypes}
The tree $\tl$ contains paths of type $\Dt$ if and only if, for every 
$\Delta_n$, $1\leq n\leq t$, every $\omega_i$ required by $\Delta_n$ 
is provided by some $\Delta_k$ with $k<n$.
\end{lemma}

The ``only if'' part of the equivalence is immediate from the preceding
discussion: the sequence $\mu_0,\mu_1,\ldots$ starts at $\mu_0=0$, and if
it moves towards any wall of the Weyl chamber before first moving 
away from it, it will pass through the wall and some $\mu_i$ will be
outside the chamber.  Conversely, if $\Dt$ is any path-type which
satisfies the condition of the lemma, then $\Dtmt$ will definitely appear
in the tree when $m_1 \gg m_2 \gg\cdots\gg m_t$.  This ensures that the
coordinates of the $\mu_i$ are always nonnegative, since the sequence of
$\mu$s moves sufficiently far away from any wall of the Weyl chamber
before the first time it moves back towards it.  We could compute the
exact conditions on the $m_i$ for a specific path; in general, they all
require that $m_n$ be bounded by some linear combination of
$m_1,\ldots,m_{n-1}$, and the first $m_i$ appearing with nonzero
coefficient in that linear combination has positive coefficient.

Now we can show that the number of nodes of path-type $\Dt$ appearing on
the $m$th level of the tree grows as $m^{t-1}$.  Consider the path
$\Dtmt$ as a point $(m_1,\ldots,m_t)$ in $\R^t$.  The path ends on row
$m$ if $m=m_1+\cdots+m_t$, so solutions lie on a plane of dimension
$t-1$; the number of solutions to that equality in nonnegative integers
is ${m+t-1}\choose{t-1}$, which certainly grows as $m^{t-1}$, as
expected.  The further linear inequalities on the $m_i$ which ensure that
$\mu_1,\ldots,\mu_m$ remain in the Weyl chamber correspond to hyperplanes
through the origin which our solutions must lie on one side of, but the
resulting region still has full dimension $t-1$ since the generic point
with $m_1 \gg m_2 \gg\cdots\gg m_t$ satisfies all of the inequalities, as
shown above.

% Inserted paragraph is here:

The highest weight of the $\g$-module at the node associated with the
generic solution of the form $m_1 \gg m_2 \gg\cdots\gg m_t$ grows
linearly in $m$.  Its dimension, therefore, grows as a polynomial in $m$,
and the degree of the polynomial is just the number of positive roots of
the Lie algebra which are not orthogonal to the highest weight.  The only
positive roots perpendicular to this generic highest weight are those
perpendicular to every highest weight which comes from a path of type
$\Dt$, and the number of such roots is the degree of polynomial growth of
the dimensions of the representations of the $\g$-module.  This can also
be expressed as $\frac12\dim(\O_{\lambda})$, where $\lambda$ is a weight
orthogonal to any given positive root if and only if all highest weights
of type $\Dt$ are.

% End inserted paragraph.

We can figure out how the multiplicities of nodes with a specific 
path-type grow as well.  Theorem~\ref{thm_decomp} gives a formula for 
multiplicities as a product of binomial coefficients over $1\leq k\leq 
r$ and $n\geq 1$.  The only terms in the product which are not 1 
correspond to nonzero values of $\delta_n - \delta_{n+1}$.  In the 
path $\Dtmt$, these occur only when $n=m_1+\cdots+m_i$ for some $1\leq 
i\leq t$, so that $\delta_n - \delta_{n+1}$ is $\Delta_i - 
\Delta_{i+1}$ (where $\Delta_{t+1}$ is just 0).  Following our 
previous notation, let $\delta_n - \delta_{n+1} = \d_n = \sum 
\d^{(k)}_n \alpha_k$.  If we take any $k$ for which $\d^{(k)}_n$ is 
nonzero, there are two possibilities for the contribution to the 
multiplicity from its binomial coefficient.  If $\omega_k$ has been 
provided by at least one of $\Delta_1,\ldots,\Delta_i$, then the value 
$P^{(k)}_n$ is a linear combination of $m_1,\ldots,m_i$, which grows 
linearly as $m$ gets large.  In this case, the binomial coefficient 
grows as a polynomial in $m$ of degree $\d^{(k)}_n$.  On the other 
hand, if $\omega_k$ has not been provided, then the binomial 
coefficient is just 1.

For any $k$, $1\leq k\leq r$, define $f(k)$ to be the smallest $i$ in 
our path-type such that $\Delta_i$ provides $\omega_k$; we say that 
$\Delta_i$ provides $\omega_k$ for the first time.  Then the total 
contribution to the multiplicity from the coordinate $k$ will be the 
product of the contributions when $n=m_1+\cdots+m_j$ for $j=f(k), 
f(k)+1, \ldots, t$.  As $m$ gets large, the product of these 
contributions grows as a polynomial of degree $\sum_{j=f(k)}^t 
\d^{(m_1+\cdots+m_j)}_n$; that is, the sum of the decreases in the 
$\alpha_k$-coordinate of the $\Delta$s.  But since $\Delta_{t+1}$ is 
just 0, that sum is exactly the $\alpha_k$-coordinate of 
$\Delta_{f(k)}$.

So given a path-type $\Dt$ which Lemma~\ref{lem_pathtypes} says appears
in $\tl$, the total of the multiplicities of the nodes of that path-type
which appear in the top $m$ rows of $\tl$ grows as a polynomial of degree
\begin{equation}
\label{defg}
g(\Dt) = t + \sum_{k=1}^{r} \alpha_k\mbox{-coordinate of }\Delta_{f(k)}
\end{equation}
where we take $\Delta_{f(k)}$ to be 0 if $\omega_k$ is not provided by any
$\Delta$ in the path-type.  This value is just the sum of the degrees of
the polynomial growths described above.

Finally, since there are only finitely many path-types, the growth of the
entire tree $\tl$ is the same as the growth of the part corresponding to
any path-type $\Dt$ which maximizes $g(\Dt)$.  So we have proven the
following, up to some calculation:

\begin{thm}
\label{thm_growth}
Let $\g$ be simply-laced with decompositions of $\wml$ given by
Theorem~\ref{thm_decomp}.  Then the dimension of the representation
$\wml$ as $m$ gets large is asymptotic to a polynomial in $m$ of degree
$\frac12\dim(\O_{\lambda})+g(\Dt)$, where the path-type $\Dt$ is one
which maximizes the value of $g$, and $\O_{\lambda}$ is the adjoint
orbit of a weight $\lambda$ which is orthogonal to exactly those positive
roots orthogonal to all highest weights of nodes with path-type $\Dt$.

\begin{enumerate}
\item
If $\g$ is of type $A_n$ then the maximum value of $g(\Dt)$ is $0$, for 
all $1\leq\l\leq n$.
\item
If $\g$ is of type $D_n$ then the maximum value of $g(\Dt)$ is 
$\lfloor\l/2\rfloor$, for $1\leq\l\leq n-2$, and $0$ for $\l=n-1,n$.
\item
If $\g$ is of type $E_6$, $E_7$, or $E_8$, the maximum value of 
$g(\Dt)$ is
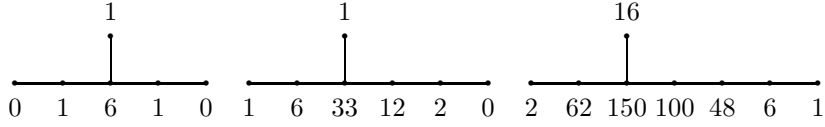
\begin{center} \setlength{\unitlength}{.25in}
\newcommand{\num}[1]{{\makebox(0,1){${#1}$}}} % bigger labels for 10pt
\newcommand{\vx}{{\circle*{.135}}}
\begin{picture}(4,3) % E_6
\multiput(0,1)(1,0){5}{\vx}\put(2,2){\vx}
\put(0,1){\line(1,0){4}}\put(2,1){\line(0,1){1}}
\put(2,2){\num{1}}\put(0,0){\num{0}}\put(1,0){\num{1}}
\put(2,0){\num{6}}\put(3,0){\num{1}}\put(4,0){\num{0}}
\end{picture}\hfill
\begin{picture}(5,3) % E_7
\put(2,2){\vx}\multiput(0,1)(1,0){6}{\vx}
\put(0,1){\line(1,0){5}}\put(2,1){\line(0,1){1}}
\put(2,2){\num{1}}\put(0,0){\num{1}}\put(1,0){\num{6}}\put(2,0){\num{33}}
\put(3,0){\num{12}}\put(4,0){\num{2}}\put(5,0){\num{0}}
\end{picture}\hfill
\begin{picture}(6,3) % E_8
\put(2,2){\vx}\multiput(0,1)(1,0){7}{\vx}
\put(0,1){\line(1,0){6}}\put(2,1){\line(0,1){1}}
\put(2,2){\num{16}}\put(0,0){\num{2}}\put(1,0){\num{62}}\put(2,0){\num{150}}
\put(3,0){\num{100}}\put(4,0){\num{48}}\put(5,0){\num{6}}\put(6,0){\num{1}}
\end{picture}\hspace*{\fill}
\end{center}
\end{enumerate}
\end{thm}

We will complete the proof by exhibiting the path-types which give 
the indicated values of $g$ and proving they are maximal.

If $\Dt$ maximizes the value of $g$, then it cannot be obtained from any
other path-type by insterting an extra $\Delta$, since any insertion
would increase the length $t$ and would not decrease the sum in the
definition of $g$.  Therefore each $\Delta_k$ in our desired path-type
must be in the positive root lattice, allowable according to
Lemma~\ref{lem_pathtypes}, and must be maximal (under $\preceq$) in
meeting those requirements; we will call a path-type maximal if this is
the case.

In particular, if $\wl$ is in the root lattice then $\Delta_1$ will
be $\wl$, and a $g$-value of 0 corresponds exactly to an $\wl$ which is
not in the root lattice and is a minimal weight.  Thus the 0s above
can be verified by inspection; these are exactly the cases in which
$\wml$ remains irreducible as a $\g$-module.  Similarly, if $\wl$ is not
in the root lattice but there is only one point in the lattice and in
the Weyl chamber under $\wl$, the path-type will consist just of that
point.  We can now limit ourselves to path-types of length greater than
one.

If $\g$ is of type $D_n$ then for each $\wl$, $2\leq\l\leq n-2$, 
there is a unique maximal path-type:
$$
\begin{array}{ll}
\wl \succ \wl-\omega_2 \succ \wl-\omega_4 \succ \cdots \succ \wl-\omega_{\l-2}
  & \mbox{when $\l$ is even} \\
\wl-\omega_1 \succ \wl-\omega_3 \succ \cdots \succ \wl-\omega_{\l-2}
  & \mbox{when $\l$ is odd}
\end{array}
$$
In both cases, the only contribution to $g$ comes from the length 
of the path, which is $\lfloor\l/2\rfloor$.  This also means that the 
nodes of the tree $\tl$ will all have multiplicity 1 in this case.

When $\g$ is of type $E_6$, $E_7$ or $E_8$, the following weights have a
unique maximal path-type (of length $>1$), whose $g$-value is given in
Theorem~\ref{thm_growth}:
$$
\begin{array}{lll}
E_6
& \l=4 & \omega_4 \succ \omega_4-\omega_2 \succ \omega_4-\omega_1-\omega_6
       \succ \omega_2+\omega_4-\omega_3-\omega_5 \succ 2\omega_2-\omega_4 \\
E_7
& \l=3 & \omega_3 \succ \omega_3-\omega_1 \succ \omega_3-\omega_6 \succ
       \omega_1+\omega_6-\omega_4 \succ 2\omega_1-\omega_3 \\
& \l=6 & \omega_6 \succ \omega_6-\omega_1 \\
E_8
& \l=1 & \omega_1 \succ \omega_1-\omega_8 \\
& \l=7 & \omega_7 \succ \omega_7-\omega_8 \succ \omega_7-\omega_1 \succ
       \omega_7+\omega_8-\omega_6 \succ 2\omega_8-\omega_7\\
& \l=8 & \omega_8
\end{array}
$$

We will consider the remaining weights in $E_8$ next.  Consider the
incomplete path-type
$$
\wl \succ \wl-\omega_8 \succ \wl-\omega_1 \succ \wl-\omega_6+\omega_8
\succ \wl+\omega_1-\omega_4+\omega_8 \succ \cdots
$$
where $\wl$ is any fundamental weight which is in the root lattice 
and high enough that all of the weights in question lie in the Weyl 
chamber.  The path so far provides $\omega_8$, $\omega_1$, $\omega_6$ 
and $\omega_4$; notice that for any $\omega_i$ which has not been 
provided, all of its neighbors in the Dynkin diagram have.  Therefore 
we can extend this path four more steps by subtracting one of 
$\alpha_2$, $\alpha_3$, $\alpha_5$ and $\alpha_7$ at each step, to 
produce a path in which every $\omega_i$ has been provided.  This can 
be extended to a full path-type by subtracting any $\alpha_i$ at each 
stage until we reach the walls of the Weyl chamber.

The resulting path-type is maximal, and is the unique maximal one up to a
sequence of transformations of the form
$$
\cdots \succ \Delta \succ \Delta-\lambda \succ \Delta-\lambda-\mu \succ \cdots
\mapsto
\cdots \succ \Delta \succ \Delta-\mu \succ \Delta-\lambda-\mu \succ \cdots
$$
which do not affect the rate of growth $g$.  All relevant weights are in
the Weyl chamber if and only if $\wl\succ\xi=(4,8,10,14,12,8,6,2)$; this
turns out to be everything except $\omega_1$, $\omega_7$ and $\omega_8$,
whose path-types are given above.  If the path-type could start at $\xi$,
it would have growth $g=8$, though this is not possible since the last
weight in the path-type would be $0$ in this case.  But each increase of
the starting point of the path by any $\alpha_i$ increases $g$ by 2
(1 from the length of the path and 1 from the multiplicity).  So the growth
for any $\wl\succ\xi$ is a linear function of its height with coefficient
2; $g=2\hyt(\wl)-120$.

The only remaining cases are $\omega_4$ and $\omega_5$ when $\g$ is of
type $E_7$.  Both work like the general case for $E_8$, beginning instead
with the incomplete path-types
$$
\begin{array}{ll}
\l=4 & \omega_4 \succ \omega_4-\omega_1 \succ \omega_4-\omega_6
     \succ \omega_1 \succ \cdots \\
\l=5 & \omega_5-\omega_7 \succ \omega_5-\omega_2 \succ
     \omega_5+\omega_7-\omega_1-\omega_6 \succ
     \omega_2+\omega_7-\omega_3 \succ \cdots
\end{array}
$$
This concludes the proof of Theorem~\ref{thm_growth}.

The same argument used for $E_8$ shows that for any choice of $\g$, all
``sufficiently large'' weights $\wl$ in a particular translate of the
root lattice will have growth given by $2\hyt(\wl)-c$ for some fixed
$c$.  A weight is sufficiently large if every $\omega_i$ is provided in
its maximal path.  Thus we can easily check that $\omega_4$ and $\omega_5$
qualify for $E_7$, and in both cases $c=63$.  Similarly, $\omega_4$ for
$E_6$ qualifies, and $c=36$.  While there are no sufficiently large
fundamental weights for $A_n$ or $D_n$, we can compute what the maximal
path-type would be if one did exist, and in all cases, $c$ is the number
of positive roots.  A uniform explanation of this fact would be nice,
even though the exhaustive computation does provide a complete proof.

\section{Computations}
\label{sec_table}

This section gives the decompositions of $\wml$ into $\g$-modules
predicted by the conjectural formulas in~\cite{KR}.  We also give the
tree structure defined in Section~\ref{sec_algorithm}.

The representations $\wml$ when $m=1$ are called fundamental
representations.  In the setting of $U_q(\g)$-module decompositions of
$U_q(\ghat)$ modules, the decompositions of the fundamental
representations for all $\g$ and most choices of $\wl$ appear
in~\cite{ChP}, calculated using techniques unrelated to the conjecture
used in~\cite{KR} to give formulas~(\ref{defZ}) and~(\ref{defPgen}).
Those computations agree with the ones given below.  In particular, the
choices of $\wl$ not calculated in~\cite{ChP} are exactly those in which
the maximal path-type (Theorem~\ref{thm_growth}) is not unique.

\subsection{\boldmath $A_n$}

As already noted, when $\g$ is of type $A_n$, the $Y(\g)$-modules $\wml$
remain irreducible when viewed as $\g$-modules.

\subsection{\boldmath $D_n$}

Let $\g$ be of type $D_n$.  As already noted, the fundamental weights
$\omega_{n-1}$ and $\omega_{n}$ are minimal with respect to $\preceq$,
so $W_m(n-1)$ and $W_m(n)$ remain irreducible as $\g$-modules.
Now suppose $\l\leq n-2$.  Then the structure of the weights in the
Weyl chamber under $\omega_\l$ does not depend on $n$, and so the
decomposition of $\wml$ in $D_n$ is the same for any $n \geq \l+2$.

As mentioned in the proof of Theorem~\ref{thm_growth}, there is a
unique maximal path-type for each $\wl$, and there are no multiplicities
greater than 1.  The decomposition is therefore very simple: if
$\l\leq n-2$ is even, then
$$
\wml \iso \Dsum_{k_2+k_4+\ldots+k_{\l-2}+k_\l = k \leq m}
     V_{ k_2 \omega_2 + k_4 \omega_4 + \ldots + k_{\l-2} \omega_{\l-2}
         + (m-k)\omega_\l }
$$
and if $\l$ is odd, then
$$
\wml \iso \Dsum_{k_1+k_3+\ldots+k_{\l-2} = k \leq m}
     V_{ k_1 \omega_1 + k_3 \omega_3 + \ldots + k_{\l-2} \omega_{\l-2}
         + (m-k)\omega_\l }
$$
where the minor difference is because $\omega_\l$ for $\l$ odd is not in
the root lattice.  The sum $k$ is the level of the tree on which that
module appears, and the parent of a module is obtained by subtracting 
1 from the first of $k_{\l-2}, k_{\l-4},\ldots$ which is nonzero (or from
$k_\l$ if nothing else is nonzero and $\l$ is even).

\subsection{\boldmath $E_n$}

When $\g$ is of type $E_n$ the tree structure is much more irregular:
these are the only cases in which a $\g$-module can appear in more that
one place in the tree and in which a node on the tree can have
multiplicitiy greater than one.

We indicate the tree structure as follows: we list every node in the
tree, starting with the root and in depth-first order, and a node on
level $k$ of the tree is written as $\oplustag{k} V_\lambda$.  This is
enough information to recover the entire tree, since the parent of
that node is the most recent summand of the form $\oplustag{k-1} V_\mu$.
Comparing Figure~\ref{fig_tree} to its representation here should make
the notation clear.

Due to space considerations, for $E_6$ we list calculations for $m\leq 
3$, for $E_7$ we list $m\leq 2$, and for $E_8$ only $m=1$.  The 
tree decomposition for $W_3(4)$ for $E_7$, for example, would have 
836 components.

\begin{enumerate} \sloppy
\item[{\boldmath $E_6$}]
\item[$W_m(1)$] remains irreducible for all $m$.
\item[$W_1(2)$] $\iso V_{\omega_2}
 \oplustag{1} V_{0}$
\item[$W_2(2)$] $\iso V_{2 \omega_2}
 \oplustag{1} V_{\omega_2}
 \oplustag{2} V_{0}$
\item[$W_3(2)$] $\iso V_{3 \omega_2}
 \oplustag{1} V_{2 \omega_2}
 \oplustag{2} V_{\omega_2}
 \oplustag{3} V_{0}$
\item[$W_1(3)$] $\iso V_{\omega_3}
 \oplustag{1} V_{\omega_6}$
\item[$W_2(3)$] $\iso V_{2 \omega_3}
 \oplustag{1} V_{\omega_3+\omega_6}
 \oplustag{2} V_{2 \omega_6}$
\item[$W_3(3)$] $\iso V_{3 \omega_3}
 \oplustag{1} V_{2 \omega_3+\omega_6}
 \oplustag{2} V_{\omega_3+2 \omega_6}
 \oplustag{3} V_{3 \omega_6}$
\item[$W_1(4)$] $\iso V_{\omega_4}
 \oplustag{1} V_{\omega_1+\omega_6}
 \oplustag{1} 2 V_{\omega_2}
 \oplustag{1} V_{0}$
\item[$W_2(4)$] $\iso V_{2 \omega_4}
 \oplustag{1} V_{\omega_1+\omega_4+\omega_6}
 \oplustag{2} V_{2 \omega_1+2 \omega_6}
 \oplustag{1} 2 V_{\omega_2+\omega_4}
 \oplustag{2} V_{\omega_3+\omega_5}
 \oplustag{2} 2 V_{\omega_1+\omega_2+\omega_6}
 \oplustag{2} 3 V_{2 \omega_2}
 \oplustag{3} V_{\omega_4}
 \oplustag{1} V_{\omega_4}
 \oplustag{2} V_{\omega_1+\omega_6}
 \oplustag{2} 2 V_{\omega_2}
 \oplustag{2} V_{0}$
\item[$W_3(4)$] $\iso V_{3 \omega_4}
 \oplustag{1} V_{\omega_1+2 \omega_4+\omega_6}
 \oplustag{2} V_{2 \omega_1+\omega_4+2 \omega_6}
 \oplustag{3} V_{3 \omega_1+3 \omega_6}
 \oplustag{1} 2 V_{\omega_2+2 \omega_4}
 \oplustag{2} V_{\omega_3+\omega_4+\omega_5}
 \oplustag{2} 2 V_{\omega_1+\omega_2+\omega_4+\omega_6}
 \oplustag{3} V_{\omega_1+\omega_3+\omega_5+\omega_6}
 \oplustag{3} 2 V_{2 \omega_1+\omega_2+2 \omega_6}
 \oplustag{2} 3 V_{2 \omega_2+\omega_4}
 \oplustag{3} V_{2 \omega_4}
 \oplustag{3} 2 V_{\omega_2+\omega_3+\omega_5}
 \oplustag{3} 3 V_{\omega_1+2 \omega_2+\omega_6}
 \oplustag{4} V_{\omega_1+\omega_4+\omega_6}
 \oplustag{3} 4 V_{3 \omega_2}
 \oplustag{4} 2 V_{\omega_2+\omega_4}
 \oplustag{1} V_{2 \omega_4}
 \oplustag{2} V_{\omega_1+\omega_4+\omega_6}
 \oplustag{3} V_{2 \omega_1+2 \omega_6}
 \oplustag{2} 2 V_{\omega_2+\omega_4}
 \oplustag{3} V_{\omega_3+\omega_5}
 \oplustag{3} 2 V_{\omega_1+\omega_2+\omega_6}
 \oplustag{3} 3 V_{2 \omega_2}
 \oplustag{4} V_{\omega_4}
 \oplustag{2} V_{\omega_4}
 \oplustag{3} V_{\omega_1+\omega_6}
 \oplustag{3} 2 V_{\omega_2}
 \oplustag{3} V_{0}$
\item[$W_1(5)$] $\iso V_{\omega_5}
 \oplustag{1} V_{\omega_1}$
\item[$W_2(5)$] $\iso V_{2 \omega_5}
 \oplustag{1} V_{\omega_1+\omega_5}
 \oplustag{2} V_{2 \omega_1}$
\item[$W_3(5)$] $\iso V_{3 \omega_5}
 \oplustag{1} V_{\omega_1+2 \omega_5}
 \oplustag{2} V_{2 \omega_1+\omega_5}
 \oplustag{3} V_{3 \omega_1}$
\item[$W_m(6)$] remains irreducible for all $m$.

\item[{\boldmath $E_7$}]
\item[$W_1(1)$] $\iso V_{\omega_1}
 \oplustag{1} V_{0}$
\item[$W_2(1)$] $\iso V_{2 \omega_1}
 \oplustag{1} V_{\omega_1}
 \oplustag{2} V_{0}$
\item[$W_1(2)$] $\iso V_{\omega_2}
 \oplustag{1} V_{\omega_7}$
\item[$W_2(2)$] $\iso V_{2 \omega_2}
 \oplustag{1} V_{\omega_2+\omega_7}
 \oplustag{2} V_{2 \omega_7}$
\item[$W_1(3)$] $\iso V_{\omega_3}
 \oplustag{1} V_{\omega_6}
 \oplustag{1} 2 V_{\omega_1}
 \oplustag{1} V_{0}$
\item[$W_2(3)$] $\iso V_{2 \omega_3}
 \oplustag{1} V_{\omega_3+\omega_6}
 \oplustag{2} V_{2 \omega_6}
 \oplustag{1} 2 V_{\omega_1+\omega_3}
 \oplustag{2} V_{\omega_4}
 \oplustag{2} 2 V_{\omega_1+\omega_6}
 \oplustag{2} 3 V_{2 \omega_1}
 \oplustag{3} V_{\omega_3}
 \oplustag{1} V_{\omega_3}
 \oplustag{2} V_{\omega_6}
 \oplustag{2} 2 V_{\omega_1}
 \oplustag{2} V_{0}$
\item[$W_1(4)$] $\iso V_{\omega_4}
 \oplustag{1} V_{\omega_1+\omega_6}
 \oplustag{1} 2 V_{\omega_2+\omega_7}
 \oplustag{1} V_{2 \omega_1}
 \oplustag{1} 3 V_{\omega_3}
 \oplustag{2} V_{\omega_6}
 \oplustag{1} V_{2 \omega_7}
 \oplustag{1} 3 V_{\omega_6}
 \oplustag{2} V_{\omega_1}
 \oplustag{1} 3 V_{\omega_1}
 \oplustag{2} V_{0}
 \oplustag{1} V_{0}$
\item[$W_2(4)$] $\iso V_{2 \omega_4}
 \oplustag{1} V_{\omega_1+\omega_4+\omega_6}
 \oplustag{2} V_{2 \omega_1+2 \omega_6}
 \oplustag{1} 2 V_{\omega_2+\omega_4+\omega_7}
 \oplustag{2} V_{\omega_3+\omega_5+\omega_7}
 \oplustag{2} 2 V_{\omega_1+\omega_2+\omega_6+\omega_7}
 \oplustag{2} 3 V_{2 \omega_2+2 \omega_7}
 \oplustag{3} V_{\omega_4+2 \omega_7}
 \oplustag{1} V_{2 \omega_1+\omega_4}
 \oplustag{2} V_{3 \omega_1+\omega_6}
 \oplustag{2} V_{4 \omega_1}
 \oplustag{1} 3 V_{\omega_3+\omega_4}
 \oplustag{2} 2 V_{\omega_1+\omega_2+\omega_5}
 \oplustag{2} 4 V_{\omega_1+\omega_3+\omega_6}
 \oplustag{2} V_{2 \omega_5}
 \oplustag{2} V_{2 \omega_2+\omega_6}
 \oplustag{2} 3 V_{\omega_4+\omega_6}
 \oplustag{3} V_{\omega_1+2 \omega_6}
 \oplustag{2} 2 V_{2 \omega_1+\omega_2+\omega_7}
 \oplustag{2} 6 V_{\omega_2+\omega_3+\omega_7}
 \oplustag{3} 2 V_{\omega_1+\omega_5+\omega_7}
 \oplustag{3} 2 V_{\omega_2+\omega_6+\omega_7}
 \oplustag{2} 3 V_{2 \omega_1+\omega_3}
 \oplustag{2} 6 V_{2 \omega_3}
 \oplustag{3} 3 V_{\omega_1+\omega_4}
 \oplustag{3} V_{\omega_2+\omega_5}
 \oplustag{3} V_{2 \omega_1+\omega_6}
 \oplustag{3} 3 V_{\omega_3+\omega_6}
 \oplustag{4} V_{2 \omega_6}
 \oplustag{1} V_{\omega_4+2 \omega_7}
 \oplustag{2} V_{\omega_1+\omega_6+2 \omega_7}
 \oplustag{2} 2 V_{\omega_2+3 \omega_7}
 \oplustag{2} V_{4 \omega_7}
 \oplustag{1} 3 V_{\omega_4+\omega_6}
 \oplustag{2} 2 V_{\omega_2+\omega_3+\omega_7}
 \oplustag{2} 3 V_{\omega_1+2 \omega_6}
 \oplustag{2} 4 V_{\omega_1+\omega_5+\omega_7}
 \oplustag{2} V_{2 \omega_3}
 \oplustag{2} V_{\omega_1+2 \omega_2}
 \oplustag{2} 3 V_{\omega_1+\omega_4}
 \oplustag{3} V_{2 \omega_1+\omega_6}
 \oplustag{2} 8 V_{\omega_2+\omega_6+\omega_7}
 \oplustag{3} 2 V_{\omega_3+2 \omega_7}
 \oplustag{2} 6 V_{\omega_2+\omega_5}
 \oplustag{3} 2 V_{\omega_3+\omega_6}
 \oplustag{3} 2 V_{\omega_1+\omega_2+\omega_7}
 \oplustag{2} V_{2 \omega_1+2 \omega_7}
 \oplustag{2} 3 V_{\omega_3+2 \omega_7}
 \oplustag{3} V_{\omega_6+2 \omega_7}
 \oplustag{2} 3 V_{2 \omega_1+\omega_6}
 \oplustag{3} V_{3 \omega_1}
 \oplustag{2} 9 V_{\omega_3+\omega_6}
 \oplustag{3} 4 V_{\omega_1+\omega_2+\omega_7}
 \oplustag{3} 3 V_{2 \omega_6}
 \oplustag{3} 4 V_{\omega_5+\omega_7}
 \oplustag{3} 4 V_{\omega_1+\omega_3}
 \oplustag{3} V_{2 \omega_2}
 \oplustag{3} 3 V_{\omega_4}
 \oplustag{4} V_{\omega_1+\omega_6}
 \oplustag{2} 3 V_{\omega_6+2 \omega_7}
 \oplustag{2} 6 V_{2 \omega_6}
 \oplustag{3} 3 V_{\omega_5+\omega_7}
 \oplustag{3} V_{\omega_4}
 \oplustag{3} V_{\omega_1+2 \omega_7}
 \oplustag{3} 3 V_{\omega_1+\omega_6}
 \oplustag{4} V_{2 \omega_1}
 \oplustag{1} 3 V_{\omega_1+\omega_4}
 \oplustag{2} 2 V_{\omega_2+\omega_5}
 \oplustag{2} 3 V_{2 \omega_1+\omega_6}
 \oplustag{2} 4 V_{\omega_3+\omega_6}
 \oplustag{2} 8 V_{\omega_1+\omega_2+\omega_7}
 \oplustag{3} 2 V_{\omega_5+\omega_7}
 \oplustag{2} V_{2 \omega_6}
 \oplustag{2} 3 V_{\omega_5+\omega_7}
 \oplustag{3} V_{\omega_1+2 \omega_7}
 \oplustag{2} 3 V_{3 \omega_1}
 \oplustag{2} 12 V_{\omega_1+\omega_3}
 \oplustag{3} 4 V_{\omega_4}
 \oplustag{3} 4 V_{\omega_1+\omega_6}
 \oplustag{2} 3 V_{2 \omega_2}
 \oplustag{3} V_{\omega_4}
 \oplustag{2} 5 V_{\omega_4}
 \oplustag{3} 3 V_{\omega_1+\omega_6}
 \oplustag{3} 4 V_{\omega_2+\omega_7}
 \oplustag{3} V_{2 \omega_1}
 \oplustag{3} 3 V_{\omega_3}
 \oplustag{4} V_{\omega_6}
 \oplustag{2} 3 V_{\omega_1+2 \omega_7}
 \oplustag{2} 9 V_{\omega_1+\omega_6}
 \oplustag{3} 4 V_{\omega_2+\omega_7}
 \oplustag{3} 3 V_{2 \omega_1}
 \oplustag{3} 4 V_{\omega_3}
 \oplustag{3} V_{2 \omega_7}
 \oplustag{3} 3 V_{\omega_6}
 \oplustag{4} V_{\omega_1}
 \oplustag{2} 6 V_{2 \omega_1}
 \oplustag{3} 3 V_{\omega_3}
 \oplustag{3} V_{\omega_6}
 \oplustag{3} 3 V_{\omega_1}
 \oplustag{4} V_{0}
 \oplustag{1} V_{\omega_4}
 \oplustag{2} V_{\omega_1+\omega_6}
 \oplustag{2} 2 V_{\omega_2+\omega_7}
 \oplustag{2} V_{2 \omega_1}
 \oplustag{2} 3 V_{\omega_3}
 \oplustag{3} V_{\omega_6}
 \oplustag{2} V_{2 \omega_7}
 \oplustag{2} 3 V_{\omega_6}
 \oplustag{3} V_{\omega_1}
 \oplustag{2} 3 V_{\omega_1}
 \oplustag{3} V_{0}
 \oplustag{2} V_{0}$
\item[$W_1(5)$] $\iso V_{\omega_5}
 \oplustag{1} V_{\omega_1+\omega_7}
 \oplustag{1} 2 V_{\omega_2}
 \oplustag{1} 2 V_{\omega_7}$
\item[$W_2(5)$] $\iso V_{2 \omega_5}
 \oplustag{1} V_{\omega_1+\omega_5+\omega_7}
 \oplustag{2} V_{2 \omega_1+2 \omega_7}
 \oplustag{1} 2 V_{\omega_2+\omega_5}
 \oplustag{2} V_{\omega_3+\omega_6}
 \oplustag{2} 2 V_{\omega_1+\omega_2+\omega_7}
 \oplustag{2} 3 V_{2 \omega_2}
 \oplustag{3} V_{\omega_4}
 \oplustag{1} 2 V_{\omega_5+\omega_7}
 \oplustag{2} V_{\omega_4}
 \oplustag{2} 2 V_{\omega_1+2 \omega_7}
 \oplustag{2} 2 V_{\omega_1+\omega_6}
 \oplustag{2} 4 V_{\omega_2+\omega_7}
 \oplustag{3} V_{\omega_3}
 \oplustag{2} 3 V_{2 \omega_7}
 \oplustag{3} V_{\omega_6}$
\item[$W_1(6)$] $\iso V_{\omega_6}
 \oplustag{1} V_{\omega_1}
 \oplustag{1} V_{0}$
\item[$W_2(6)$] $\iso V_{2 \omega_6}
 \oplustag{1} V_{\omega_1+\omega_6}
 \oplustag{2} V_{2 \omega_1}
 \oplustag{1} V_{\omega_6}
 \oplustag{2} V_{\omega_1}
 \oplustag{2} V_{0}$
\item[$W_m(7)$] remains irreducible for all $m$.

\item[{\boldmath $E_8$}]
\item[$W_1(1)$] $\iso V_{\omega_1}
 \oplustag{1} V_{\omega_8}
 \oplustag{1} V_{0}$
\item[$W_1(2)$] $\iso V_{\omega_2}
 \oplustag{1} V_{\omega_7}
 \oplustag{1} 2 V_{\omega_1}
 \oplustag{1} 2 V_{\omega_8}
 \oplustag{1} V_{0}$
\item[$W_1(3)$] $\iso V_{\omega_3}
 \oplustag{1} V_{\omega_6}
 \oplustag{1} 2 V_{\omega_1+\omega_8}
 \oplustag{1} 3 V_{\omega_2}
 \oplustag{2} V_{\omega_7}
 \oplustag{1} V_{2 \omega_8}
 \oplustag{1} 3 V_{\omega_7}
 \oplustag{2} V_{\omega_1}
 \oplustag{1} 4 V_{\omega_1}
 \oplustag{2} 2 V_{\omega_8}
 \oplustag{1} 3 V_{\omega_8}
 \oplustag{2} V_{0}
 \oplustag{1} V_{0}$
\item[$W_1(4)$] $\iso V_{\omega_4}
 \oplustag{1} V_{\omega_1+\omega_6}
 \oplustag{1} 2 V_{\omega_2+\omega_7}
 \oplustag{1} V_{2 \omega_1+\omega_8}
 \oplustag{1} 3 V_{\omega_3+\omega_8}
 \oplustag{2} V_{\omega_6+\omega_8}
 \oplustag{1} V_{2 \omega_7}
 \oplustag{1} 6 V_{\omega_1+\omega_2}
 \oplustag{2} 2 V_{\omega_5}
 \oplustag{2} 2 V_{\omega_1+\omega_7}
 \oplustag{1} 3 V_{\omega_6+\omega_8}
 \oplustag{2} V_{\omega_1+2 \omega_8}
 \oplustag{1} 5 V_{\omega_5}
 \oplustag{2} 3 V_{\omega_1+\omega_7}
 \oplustag{2} 4 V_{\omega_2+\omega_8}
 \oplustag{2} V_{2 \omega_1}
 \oplustag{2} 3 V_{\omega_3}
 \oplustag{3} V_{\omega_6}
 \oplustag{1} 3 V_{\omega_1+2 \omega_8}
 \oplustag{2} V_{3 \omega_8}
 \oplustag{1} 9 V_{\omega_1+\omega_7}
 \oplustag{2} 4 V_{\omega_2+\omega_8}
 \oplustag{2} 3 V_{2 \omega_1}
 \oplustag{2} 4 V_{\omega_3}
 \oplustag{2} 4 V_{\omega_7+\omega_8}
 \oplustag{2} 3 V_{\omega_6}
 \oplustag{3} V_{\omega_1+\omega_8}
 \oplustag{1} 10 V_{\omega_2+\omega_8}
 \oplustag{2} 4 V_{\omega_3}
 \oplustag{2} 6 V_{\omega_7+\omega_8}
 \oplustag{2} 6 V_{\omega_6}
 \oplustag{2} 8 V_{\omega_1+\omega_8}
 \oplustag{3} 2 V_{\omega_2}
 \oplustag{1} 6 V_{2 \omega_1}
 \oplustag{2} 3 V_{\omega_3}
 \oplustag{2} V_{\omega_6}
 \oplustag{2} 3 V_{\omega_1+\omega_8}
 \oplustag{3} V_{2 \omega_8}
 \oplustag{1} 7 V_{\omega_3}
 \oplustag{2} 5 V_{\omega_6}
 \oplustag{2} 8 V_{\omega_1+\omega_8}
 \oplustag{2} 9 V_{\omega_2}
 \oplustag{3} 3 V_{\omega_7}
 \oplustag{2} V_{2 \omega_8}
 \oplustag{2} 3 V_{\omega_7}
 \oplustag{3} V_{\omega_1}
 \oplustag{1} V_{3 \omega_8}
 \oplustag{1} 8 V_{\omega_7+\omega_8}
 \oplustag{2} 3 V_{\omega_6}
 \oplustag{2} 4 V_{\omega_1+\omega_8}
 \oplustag{2} 3 V_{2 \omega_8}
 \oplustag{3} V_{\omega_7}
 \oplustag{1} 7 V_{\omega_6}
 \oplustag{2} 5 V_{\omega_1+\omega_8}
 \oplustag{2} 8 V_{\omega_2}
 \oplustag{2} 3 V_{2 \omega_8}
 \oplustag{2} 9 V_{\omega_7}
 \oplustag{3} 3 V_{\omega_1}
 \oplustag{2} 3 V_{\omega_1}
 \oplustag{3} V_{\omega_8}
 \oplustag{1} 15 V_{\omega_1+\omega_8}
 \oplustag{2} 8 V_{\omega_2}
 \oplustag{2} 9 V_{2 \omega_8}
 \oplustag{2} 12 V_{\omega_7}
 \oplustag{2} 9 V_{\omega_1}
 \oplustag{3} 3 V_{\omega_8}
 \oplustag{2} 3 V_{\omega_8}
 \oplustag{3} V_{0}
 \oplustag{1} 8 V_{\omega_2}
 \oplustag{2} 6 V_{\omega_7}
 \oplustag{2} 10 V_{\omega_1}
 \oplustag{2} 6 V_{\omega_8}
 \oplustag{1} 6 V_{2 \omega_8}
 \oplustag{2} 3 V_{\omega_7}
 \oplustag{2} V_{\omega_1}
 \oplustag{2} 3 V_{\omega_8}
 \oplustag{3} V_{0}
 \oplustag{1} 7 V_{\omega_7}
 \oplustag{2} 5 V_{\omega_1}
 \oplustag{2} 8 V_{\omega_8}
 \oplustag{2} V_{0}
 \oplustag{1} 7 V_{\omega_1}
 \oplustag{2} 5 V_{\omega_8}
 \oplustag{2} 3 V_{0}
 \oplustag{1} 5 V_{\omega_8}
 \oplustag{2} 3 V_{0}
 \oplustag{1} V_{0}$
\item[$W_1(5)$] $\iso V_{\omega_5}
 \oplustag{1} V_{\omega_1+\omega_7}
 \oplustag{1} 2 V_{\omega_2+\omega_8}
 \oplustag{1} V_{2 \omega_1}
 \oplustag{1} 3 V_{\omega_3}
 \oplustag{2} V_{\omega_6}
 \oplustag{1} 2 V_{\omega_7+\omega_8}
 \oplustag{1} 4 V_{\omega_6}
 \oplustag{2} 2 V_{\omega_1+\omega_8}
 \oplustag{2} 2 V_{\omega_2}
 \oplustag{1} 6 V_{\omega_1+\omega_8}
 \oplustag{2} 2 V_{\omega_2}
 \oplustag{2} 2 V_{2 \omega_8}
 \oplustag{2} 2 V_{\omega_7}
 \oplustag{1} 5 V_{\omega_2}
 \oplustag{2} 3 V_{\omega_7}
 \oplustag{2} 4 V_{\omega_1}
 \oplustag{1} 3 V_{2 \omega_8}
 \oplustag{2} V_{\omega_7}
 \oplustag{1} 5 V_{\omega_7}
 \oplustag{2} 3 V_{\omega_1}
 \oplustag{2} 4 V_{\omega_8}
 \oplustag{1} 5 V_{\omega_1}
 \oplustag{2} 3 V_{\omega_8}
 \oplustag{2} V_{0}
 \oplustag{1} 4 V_{\omega_8}
 \oplustag{2} 2 V_{0}
 \oplustag{1} V_{0}$
\item[$W_1(6)$] $\iso V_{\omega_6}
 \oplustag{1} V_{\omega_1+\omega_8}
 \oplustag{1} 2 V_{\omega_2}
 \oplustag{1} V_{2 \omega_8}
 \oplustag{1} 3 V_{\omega_7}
 \oplustag{2} V_{\omega_1}
 \oplustag{1} 3 V_{\omega_1}
 \oplustag{2} V_{\omega_8}
 \oplustag{1} 3 V_{\omega_8}
 \oplustag{2} V_{0}
 \oplustag{1} V_{0}$
\item[$W_1(7)$] $\iso V_{\omega_7}
 \oplustag{1} V_{\omega_1}
 \oplustag{1} 2 V_{\omega_8}
 \oplustag{1} V_{0}$
\item[$W_1(8)$] $\iso V_{\omega_8}
 \oplustag{1} V_{0}$

\fussy
\end{enumerate}

\noindent
{\sc Department of Mathematics, University of California Berkeley,
Berkeley, CA, 94720, USA;} {\tt kleber@math.berkeley.edu}

\end{document}